\documentclass[11pt, a4paper]{article}
\usepackage{amsmath,amssymb,latexsym,graphics,epsfig}
\usepackage{hyperref}
\usepackage{color}
\usepackage{amsthm}

\newtheorem{theo}{Theorem}[section]

\newtheorem{lema}[theo]{Lemma}

\newtheorem{example}[theo]{Example}
\newtheorem{coro}[theo]{Corollary}

\def\proof{{\noindent\mbox{\boldmath $Proof.$}}\hskip 0.3truecm}
\def\endproof{\mbox{ \quad $\Box$}}

\def\A{\mbox{\boldmath $A$}}
\def\B{\mbox{\boldmath $B$}} 
\def\C{\mbox{\boldmath $C$}}
\def\Cb{\overline{C}}

\def\D{\mbox{\boldmath $D$}}

\def\G{\Gamma}

\def\I{\mbox{\boldmath $I$}}

\def\K{\mbox{\boldmath $K$}}
\def\L{\mbox{\boldmath $L$}}
\def\M{{\cal M}}

\def\Part{{\cal P}}

\def\S{\mbox{\boldmath $S$}}

\def\R{\mathbb R}

\def\d{\partial}
\def\diag{{\rm diag\,}}

\def\o{{\rm o}}

\def\ecc{\varepsilon}
\def\e{\mbox{\boldmath $e$}}

\def\db{\overline{\delta}}
\def\Cb{\overline{C}}
\def\Xb{\overline{X}}
\def\Yb{\overline{Y}}

\def\j{\mbox{\boldmath $j$}}

\def\vecnu{\mbox{\boldmath $\nu$}}
\def\vecrho{\mbox{\boldmath $\rho$}}
\def\vecmu{\mbox{\boldmath $\mu$}}
\def\veceta{\mbox{\boldmath $\eta$}}
\def\vec0{\mbox{\bf 0}}

\def\tr{\mathop{\rm tr}\nolimits}
\def\ev {\mathop{\rm ev }\nolimits}

\begin{document}
\title{Eigenvalue Interlacing \\ and Weight Parameters of Graphs
     \thanks{
             Work supported in part by the Spanish Research Council
(Comisi\'on Interministerial de Ciencia y Tecnolog\'\i a, CICYT)
     under projects TIC 94-0592 and TIC 97-0963. \newline \indent  } }

\author{M.A. Fiol
\\ {\small Departament de Matem\`atica Aplicada i Telem\`atica} \\
{\small Universitat Polit\`ecnica de Catalunya} \\ {\small Jordi
Girona 1-3 , M\`odul C3, Campus Nord }\\ {\small 08034 Barcelona,
Spain; email: {\tt fiol@mat.upc.es}}   }
\date{ }
\maketitle

%%%%%%%%%%%%%%%%%%%%%%%%%%%%%%%%%%%%%%%%%%%%%%%%%%%%%%%%%%
%\newpage
%\vskip 3cm
%
%\noindent{\bf Running Head:} \\
%Interlacing and weight parameters
%\vskip 1cm
%\noindent{\bf Key words:} \\
%Adjacency matrix; Eigenvalue interlacing;  Laplacian matrix
%\vskip 1cm
%\noindent{\bf AMS classifications:} \\
%05C50
%\vskip 1cm
%\noindent{\bf Mailing Address:} \\
%M.A. Fiol \\
%{\em Departament de Matem\`atica Aplicada i Telem\`atica}\\
%Universitat Polit\`ecnica de Catalunya \\
%Jordi Girona 1-3 \\
%M\`odul C3, Campus Nord \\
%08034 Barcelona, SPAIN.
%\vskip 1cm
%\noindent{\bf Email Address:} \\
%{\tt fiol@mat.upc.es}
%
%
%%%%%%%%%%%%%%%%%%%%%%%%%%%%%%%%%%%%%%%%%%%%%%%%%%%%%%%%%%%%%
%\newpage
\begin{abstract}
Eigenvalue interlacing  is a
versatile technique for deriving results in algebraic combinatorics. In particular, it has been
successfully used for proving a number of results about the relation between the
(adjacency matrix or Laplacian) spectrum of a graph and some of its properties. For
instance, some characterizations of regular partitions, and bounds
for some parameters, such as the independence and chromatic numbers, the
diameter, the bandwidth, etc., have been obtained.  For each parameter of a
graph involving the cardinality of some vertex sets, we can define its
corresponding weight parameter by giving some ``weights" (that is, the
entries of the positive eigenvector) to the vertices and replacing
cardinalities by square norms. The key point is that such weights
``regularize" the graph, and hence allow us to define a kind of regular
partition, called ``pseudo-regular," intended for general graphs. Here we
show how to use interlacing for proving results about some weight
parameters and pseudo-regular partitions of a  graph. For instance, generalizing a
well-known result of Lov\'asz, it is shown that the weight Shannon capacity $\Theta^*$ of a
connected graph $\G$, with $n$ vertices and (adjacency matrix) eigenvalues
$\lambda_1>\lambda_2\ge\cdots \ge \lambda_n$, satisfies
$$
\Theta\le \Theta^* \le \frac{\|\vecnu\|^2}{1-\frac{\lambda_1}{\lambda_n}}
$$
where $\Theta$ is the (standard) Shannon capacity and $\vecnu$ is the positive eigenvector
normalized to have smallest entry $1$. In the special case of regular graphs, the results
obtained have some interesting corollaries, such as an upper bound for some of the multiplicities
of the eigenvalues of a distance-regular graph. Finally, some results involving the Laplacian
spectrum are derived.
\end{abstract}

\newpage

%%%%%%%%%%%%%%Sec. 1
\section{Introduction}

As has been shown by Haemers \cite{ha80,ha95} and other authors, eigenvalue interlacing is a
powerful technique for deriving results about combinatorial structures from the spectra of their
associated matrices. A good  and quite complete survey on this topic is Haemers' paper
\cite{ha95}. In particular, this technique allows us to infer a number of properties of a graph,
such as bounds for its diameter, independence and chromatic numbers, bandwidth, etc, from
(part of) its spectrum. Before explaining the contents of this work, we will introduce some basic
terminology.

Let $\G=(V,E)$ be a (finite and simple) graph on $n:=|V|$
vertices. Throughout the paper, $\G$ will be supposed to be  non-trivial, that is
$E\not=\emptyset$. The distance between two vertices $u,v\in V$ will be denoted by
$\partial(u,v)$. Then, the distance
between two subsets $U_1,U_2\subset V$ is  $\partial(U_1,U_2):=\min_{u\in U_1,v\in
U_2}\{\partial(u,v)\}$. For a given vertex subset
$C\subset V$ and some integer $k\ge 0$, we denote by
$N_k(C)$ the set of vertices at distance
$k$ from (some vertex of) $C$. Similarly, $\G_k$ stands for the graph with vertex set $V$
where two vertices are adjacent whenever they are at distance $k$ in $\G$. Thus,
$N_0 (C)=C$, $N_1(\{u\})=\G(u)$, the set of vertices adjacent to $u$, and $\G_1=\G$.
The {\it eccentricity} of $C$, denoted $\ecc_C$, can be defined as the maximum distance of any
vertex of $\G$ from $C$. (In Coding Theory this correspond to the
``covering radius" of $C$.) The notation $\Cb$ will be used to denote the complement of the set
$C$ in $V$.

The eigenvalues of (the adjacency matrix $\A(\G)$ of) $\G$ will be denoted by
$\lambda_1\ge\lambda_2\ge\cdots\ge\lambda_n$ (including multiplicities). If $\G$ is
connected, the theorem of Perron-Frobenius assures that
$\lambda_1$ is simple, positive (in fact, it coincides with the {\it spectral radius} of $\A(\G)$),
and  with positive eigenvector. If $\G$ is not connected, the existence of such an eigenvector is
not guaranteed, unless all its connected components have the same maximum eigenvalue.
Throughout the paper,  it is supposed that the eigenvalue
$\lambda_1$  has indeed a positive  eigenvector,  denoted by
$\vecnu$, which is normalized in such a way that its minimum entry (in each
connected component of $\G$) is
$1$. For instance, if $\G$ is regular, we just have $\vecnu=\j$, the all-$1$ vector.
Usually, vectors and matrices are indexed by the vertices of $V$, so that
the above condition reads $\min_{u\in U}\{\nu_u\}=1$ for every  vertex set $U\subseteq
V$ of a connected component. When we are interested in the set of distinct eigenvalues, the
notation
$\ev \G
\equiv
\ev
\A(\G) =\{\theta_0>\theta_1>\cdots >\theta_d\}$ will be used (note that
$\theta_0=\lambda_1$ and $\theta_d=\lambda_n$).

Consider the map $\vecrho:\Part(V) \longrightarrow \R^n$ defined by
$\vecrho U:=\sum_{u\in U} \nu_u\e_u$ for any $U\neq \emptyset$, where $\e_u$ represents the
$u$-th canonical (column) vector,  and $\vecrho \emptyset:={\bf 0}$. Note that, with
$\vecrho u:=\vecrho\{u\}$, we have $\|\vecrho u\|=\nu_u$, so that we can see $\vecrho$ as a
function which assigns weights to the vertices of $\G$. In doing so we ``regularize" the graph, in
the sense that the {\it average weight degree} of each vertex $u\in V$ becomes a
constant:
\begin{equation}
\label{regularize}
\delta^*(u):=\frac{1}{\nu_u}\sum_{v\in\G(u)}\nu_v=\lambda_1.
\end{equation}
Using these weights, we can also consider the so-called ``pseudo-regular partitions" of a graph,
defined in the next section, which generalize the standard notion of regular (or equitable)
partitions.

In this context, the author
\cite{f1} introduced the notion of a ``weight parameter" of a graph, defined as follows. For each
parameter of a graph $\Gamma$, say $\xi$, defined as the maximum cardinality of a
set  $U\subset V$ satisfying a given property {\sf P}, we define the corresponding {\it
weight parameter}, denoted by $\xi^*$, as the maximum value of $\|\vecrho U\|^2$ of
a vertex set $U$ satisfying {\sf P}. Note that, when the graph is regular, we have $\vecnu=\j$
and then $\xi^*\equiv \xi$. Otherwise,  when we are dealing with  non-regular graphs,  the
weight parameters are sometimes  more convenient to work with, as we will see later.  As an
instance of weight parameter,  let us consider the {\it weight independence number} of
$\G$, defined as
$$
\alpha^* := \max_{U\subset V}\{\|\vecrho U\|^2 : \mbox{$U$  is an independent set}\}.
$$

The main concern of this paper is the use of eigenvalue interlacing  for obtaining results on
some weight parameters and pseudo-regular partitions. It is shown that this approach leads
sometimes to simple proofs for some results concerning standard parameters, such as the
chromatic index and the Shannon capacity of a (not necessarily regular) graph. The basic tools for
our study are explained in the following section. The remaining sections are devoted to applying
the technique in different situations where either the adjacency matrix or the  Laplacian matrix
is considered.

%%%%%%%%%%%%%%%%%%%%%%%%%%%%%%%%Sec. 2
\section{Interlacing and Pseudo-Regular Partitions}

Our starting point is the following theorem, proved by Haemers in \cite{ha80,ha95}. That
author alludes to the first part of the theorem as a classical result, referring the reader to
Courant and Hilbert's book \cite{CH}.

\begin{theo}
\label{interlac-theo}
Let $\A$ be a real symmetric $n\times n$ matrix with eigenvalues $\lambda_1\ge
\lambda_2\ge \cdots\ge \lambda_n$. For some integer $m<n$,  let $\S$ be a  real $n\times m$
matrix with orthonormal columns, that is $\S^\top\S=\I$, and consider the matrix
$\B:=\S^\top\A\S$, with eigenvalues $\mu_1\ge \mu_2\ge\cdots \ge \mu_m$. Then the
following statements hold.
\begin{itemize}
\item[(a)]
The eigenvalues of $\B$ interlace the eigenvalues of $\A$. That is,
$$
\lambda_i\ge \mu_i\ge \lambda_{n-m+i}\ \ \ (1\le i\le m).
$$
\item[(b)]
If the interlacing is tight, that is, for some  $0\le k\le m$,
$\lambda_i= \mu_i$ $(1\le i\le k)$ and $\mu_i= \lambda_{n-m+i}$  $(k+1\le i\le m)$,
then  $\S\B=\A\S$. \endproof
\end{itemize}
\end{theo}

%%%%%%%%%%%%%\subsection{Pseudo-regular partitions}
Let $\G=(V,E)$ be a graph with adjacency matrix $\A:=\A(\G)$ and positive eigenvector
$\vecnu$ with elements indexed by the vertices of $\G$. A partition $\Part$ of the vertex set
$V=V_1\cup\cdots\cup V_m$ is called {\it pseudo-regular} (or {\it pseudo-equitable})
whenever the {\it $($pseudo-$)$intersection numbers}
\begin{equation}
\label{intersec-num}
b_{ij}^*(u):=\frac{1}{\nu_u}\sum_{v\in\G(u)\cap V_j} \nu_v\ \ \ (1\le i,j\le m)
\end{equation}
 do not depend on the chosen vertex $u\in V_i$, but only on the subsets $V_i$ and $V_j$.
In this case, such numbers are simply written as $b_{ij}^*$, and the $m\times m$ matrix
$\B^*=(b_{ij}^*)$ is  referred to as  the {\it pseudo-quotient matrix} of $\A$ with respect to the
(pseudo-regular) partition $\Part$.  Pseudo-regular partitions were introduced by Garriga and
the author  \cite{fg3}, as a generalization of the so-called regular partitions, where the above
numbers are just defined by $b_{ij}^*(u):=|\G(u)\cap V_j|$ $(u\in V_i)$. A detailed study of
regular partitions can be found in Godsil \cite{g93} and Godsil and  McKay \cite{gm80}.
 (See also Brouwer, Cohen, and Neumaier \cite{bcn89} and McKay \cite{m76}.)
 A vertex subset $C\subset V$ is said
to be a {\it completely pseudo-regular code} if the distance partition around $C$, that is
$V=C\cup N_1(C)\cup\cdots\cup N_{\varepsilon_C}(C)$, is pseudo-regular.  A spectral
characterization of such codes can be found in \cite{fg3}.

Of course we can also define, in the same way, a pseudo-regular partition of (the rows and
columns of) any matrix $\A$ with a positive eigenvector. For instance, by the Perron-Frobenius
theorem, this is the case when $\A$ is an $n\times n$ non-negative irreducible matrix. (In this
case the corresponding eigenvalue is simple, non-zero ($n>1$), and coincides with the spectral
radius.) Another example is when $\A$ is the {\it Laplacian matrix} of a graph $\G$, denoted by
$\L\equiv \L(\G)$, and defined as $\L(\G):=\D-\A(\G)$, where $\D={\rm
diag}(\delta_1,\delta_2,\ldots,\delta_n)$ and $\delta_i$ stands for the degree of the $i$-th
vertex. (See  Mohar \cite{m88} for a comprehensive survey on the properties and applications of
such a matrix, and Rodr\'\i guez \cite{r97} for some recent results involving it.) Indeed, $\L$ has
eigenvalues $\lambda_1=0\le \lambda_2\le \cdots\le \lambda_n$, (they are usually enumerated
in non-decreasing order) and the  eigenvalue $0$ has the eigenvector $\j$.
Then, although in this paper we limit ourselves to the adjacency and Laplacian matrices,
most of the results obtained remain valid for ``appropriate" matrices which satisfy the
hypotheses of Theorem \ref{interlac-theo} and have a positive eigenvector.
%associated to their maximum eigenvalue.
Here ``appropriate" means that the considered matrices give some information about the
structure of the graph, which is relevant to the parameter(s) under consideration. (In other
words, matrices with an appropriate underlying graph.) For instance,  in the study of the weight
independence number, undertaken in the next section,
%%we can use any
%%As commented in Section 2,  the above results remain valid for any matrix
%%$\A^*\in \mho(\G)$, where $\mho(\G)$ is the set of (non-zero) irreducible
%%symmetric  matrix $\A$ with entries
%%$(\A)_{uv}\ge 0$ if $u\sim v$, and $(\A)_{uv}= 0$ otherwise (as in the proof of the main
%%theorem we use the implication $u\not\sim v\Rightarrow (\A)_{uv}= 0$).
we need to use a matrix $\A$ such that if $u,v$ are non-adjacent vertices then $(\A)_{uv}= 0$.

A matrix characterization of pseudo-regular partitions can be done via the following
matrix associated with (any) partition $\Part$: $V_1\cup\cdots\cup V_m$. The {\it
weight-characteristic} matrix of $\Part$ is the $n\times m$ matrix $\S^*=(s_{uj}^*)$ with
entries
$$
s_{uj}^* =\left\{ \begin{array}{ll}
\nu_u & \mbox{\rm if } u\in V_j, \\
0 & \mbox{\rm otherwise}.
\end{array}\right.
$$
\begin{lema}
Let $\G=(V,E)$ be a graph with adjacency matrix $\A$ and positive eigenvector
$\vecnu$, and consider a vertex partition $\Part$
with weight-characteristic matrix $\S^*$. Then $\Part$ is pseudo-regular if and only if there
exists an ($m\times m$) matrix $\C$ such that $\S^*\C=\A\S^*$. Moreover, in this
case  $\C=\B^*$, the pseudo-quotient matrix of $\A$ with respect to $\Part$.
\end{lema}

\proof
Let $\C=(c_{ij})$ be an $m\times m$ matrix. Let  $u\in V_i$ and $1\le j\le m$. Then, the
result follows from the equalities:

$(\S^*\C)_{uj} = \sum_{k=1}^m s_{uk}^*c_{kj}=\nu_u c_{ij}$;

$(\A\S^*)_{uj}  =  \sum_{v\in V} a_{uv} s_{vj}^*=\sum_{v\in\G(u)\cap V_j} \nu_v= \nu_u b_{ij}^*(u)$,

\noindent where we have used the definition of $b_{ij}^*(u)$.
\endproof

Most of the results about regular partitions can be generalized for pseudo-regular partitions. For
instance, using the above lemma it can be proved that all the eigenvalues of the
pseudo-quotient matrix $\B^*$ are also eigenvalues of $\A$ (see Garriga's thesis \cite{g97}).

 Let us now consider a new $n\times m$ matrix, $\S=(s_{uj})$, obtained by just normalizing the
columns of $\S^*$. Namely,
$$
s_{uj} =\left\{ \begin{array}{ll}
\frac{\nu_u}{\|\vecrho V_j\|} & \mbox{\rm if } u\in V_j, \\
0 & \mbox{\rm otherwise}
\end{array}\right.
$$
and, hence, satisfying $\S^\top\S=\I$. From such a matrix we define the {\it
weight-quotient} matrix of $\A$, with respect to $\Part$, as $\B:=\S^\top \A\S$. Notice that
this matrix  has entries
$$
b_{ij}=\sum_{u,v\in V} s_{ui}a_{uv} s_{vj}=\sum_{u\in V_i, v\in V_j}
a_{uv}\frac{\nu_u}{\|\vecrho V_i\|}\frac{\nu_v}{\|\vecrho V_j\|}
=\frac{1}{\|\vecrho V_i\|\|\vecrho V_j\|}\sum_{uv\in E(V_i, V_j)}\nu_u\nu_v=b_{ji}
$$
where $E(V_i, V_j)$ stands for the set of edges with endpoints in $V_i$ and $V_j$ (when
$V_i=V_j$ each edge counts twice). In particular, note that when $\G$ is regular
$b_{ij}=|E(V_i, V_j)|/\sqrt{|V_i||V_j|}$, so that if $|V_i|=|V_j|$ for any $1\le i,j\le m$,
then $\B$ coincides with the ``quotient matrix" used by Haemers \cite{ha95} (with entries being
the average row sums of the submatrices induced by the partition).  In the case
$N_1(V_i)\subset V_j$ we get, by (\ref{regularize}),
$$
b_{ij}
=\frac{1}{\|\vecrho V_i\|\|\vecrho V_j\|}\sum_{u\in V_i}\nu_u\sum_{v\in \G(u)}\nu_v
=\frac{\lambda_1}{\|\vecrho V_i\|\|\vecrho V_j\|}\sum_{u\in V_i}\nu_u^2
=\frac{\lambda_1\|\vecrho V_i\|}{\|\vecrho V_j\|}.
$$
In addition, we will also use the fact that  $\B$ has eigenvalue $\lambda_1$, with
corresponding eigenvector
$\vecmu:=\S^\top\vecnu=\left(\|\vecrho V_1\|, \ldots, \|\vecrho V_m\|\right)^\top$. Indeed,
\begin{eqnarray*}
(\B\vecmu)_i & = &\sum_{j=1}^m \sum_{u\in V_i, v\in V_j}
\frac{a_{uv}\nu_u\nu_v}{\|\vecrho V_i\|\|\vecrho V_j\|}\|\vecrho V_j \|
=\frac{1}{\|\vecrho V_i\|}\sum_{u\in V_i}\nu_u\sum_{v\in \G(u)}\nu_v \\
 & = & \frac{\lambda_1}{\|\vecrho V_i\|}\sum_{u\in V_i}\nu_u^2
=\lambda_1\vecmu_i \ \ \ (1\le i\le m).
\end{eqnarray*}

The following result, which is basic to our study,  is a direct consequence of Theorem
\ref{interlac-theo}, and can be thought of as a generalization of Corollary 2.3 in \cite{ha95}.

\begin{lema}
\label{interlac-coro}
Let  $\G=(V,E)$ be graph with adjacency matrix $\A$ and positive eigenvector
$\vecnu$, and consider a partition $\Part$ of $V$ inducing the weight-quotient matrix $\B$. Then
the following hold:
\begin{itemize}
\item[(a)]
The eigenvalues of $\B$ interlace the eigenvalues of $\A$;
\item[(b)]
If the interlacing is tight,  then  the partition $\Part$ is pseudo-regular.
\end{itemize}
\end{lema}

\proof
We only need to prove (b). If the interlacing is tight we already know, by Theorem
\ref{interlac-theo}(b),  that $\S\B=\A\S$. Moreover, $\S=\S^*\D$, with $\D:=\diag (\|\vecrho
V_1\|^{-1},\ldots,\|\vecrho V_m\|^{-1})$.
%=((\S^*)^\top\S^*)^{-1/2}$.
Hence,
$$
\S\S^\top\A\S=\S^*(\D\S^\top\A\S^*)\D=\A\S^*\D \ \ \Rightarrow
\ \ \S^*\B^*=\A\S^*
$$
with $\B^*:=\D\S^\top\A\S^*=\D\B\D^{-1}$
%\D\D^\top(\S^*)^\top\A\S^*$
being the pseudo-quotient matrix of
$\A$ with respect to $\Part$.
\endproof

%%%%%%%%%%%%%%%%%%%%%%%%%Sec. 3
\section{The weight independence number}

Using the results derived above, mainly Lemma \ref{interlac-coro}, most of the results
obtained for regular graphs can be extended to general graphs (with a positive eigenvector). The
only difference is that we must now consider weight parameters and pseudo-equitable
partitions. Inspired by Haemers' paper \cite{ha95},  we first derive an upper bound for both the
weight independence number and the Shannon capacity of a graph. As a straightforward
consequence of the former, we then obtain the well-known Hoffman's upper bound for the
chromatic number.

\begin{theo}
\label{theo-weight-indepen}
Let $\G$ be a graph with  eigenvalues $\lambda_1\ge\cdots \ge\lambda_n$ and
positive eigenvector $\vecnu$. Then, its weight independence number  satisfies
\begin{equation}
\label{weight-indepen}
\alpha^*\le \frac{\|\vecnu\|^2}{1-\frac{\lambda_1}{\lambda_n}}.
\end{equation}
If the bound is attained for some independent set $C$, then $C$ is a completely pseudo-regular
code with eccentricity $\ecc_C=2$.
\end{theo}

\proof
Let $C\subset V$ such that $\alpha^*=\|\vecrho C\|^2$, and let $\Part$ be the
partition $V_1\cup V_2=C\cup \Cb$, where $\Cb:=V\setminus C$. Then,
the weight-quotient matrix of $\A:=\A(\G)$ with respect to $\Part$ turns out to be
\begin{equation}
\label{B}
\B=\lambda_1\left(\begin{array}{cc}
0 & \frac{\|\vecrho C\|^2} {\|\vecrho C\|\|\vecrho\Cb\|} \\
\frac{\|\vecrho C\|^2} {\|\vecrho C\|\|\vecrho\Cb\|} &
\frac{\|\vecrho\Cb\|^2-\|\vecrho C\|^2}{\|\vecrho\Cb\|^2}
\end{array}\right)
\end{equation}
with eigenvalues $\mu_1=\lambda_1$  and
$$
\mu_2 = \tr \B - \lambda_1
=\frac{-\lambda_1\|\vecrho C\|^2}{\|\vecnu\|^2-\|\vecrho C\|^2}
=\frac{-\lambda_1\alpha^*}{\|\vecnu\|^2-\alpha^*}.
$$
Hence, since $\mu_2\ge \lambda_n$ by Lemma \ref{interlac-coro}, the result follows.
In addition, if equality holds, then the interlacing is tight (since $\mu_1=\lambda_1$
and $\mu_2=\lambda_n$) and therefore the partition is pseudo-regular. In particular, from the
corresponding pseudo-quotient matrix $\B^*=\D\B\D^{-1}$, we get that, for every vertex $u\in
\Cb$,
$$
b_{21}^*(u)=\frac{1}{\nu_u}\sum_{v\in\G(u)\cap C} \nu_v
=\frac{\lambda_1\|\vecrho C\|^2}{\|\vecrho \Cb\|^2}
=-\lambda_n\neq 0.
$$
Consequently, $\ecc_C=2$ and $\Part$ is the distance partition around $C$.
\endproof

Let $\nu_{\max}:=\max_{u\in V}\{\nu_u\}$. Then, since clearly $\alpha^*\ge \nu_{\max}^2$,
the above theorem gives
$$
1-\frac{\lambda_1}{\lambda_n}\le \frac{\|\vecnu\|^2}{\nu_{\max}^2}\le n
$$
for any such graph $\G$, with equality holding in both iff $\G$ is the complete graph $K_n$.
%We end this section by noting that, from Theorem  \ref{theo-weight-indepen}, we can conclude
%that the complete graphs $K_n$ are the only graphs satisfying
%$\|\vecnu\|^2(=n)=1-\frac{\lambda_1}{\lambda_n}$.

As another simple corollary of Theorem \ref{theo-weight-indepen} we can get the known result
of Hoffman \cite{hof70}, which provides a lower bound on  the chromatic number  $\chi$ of any
graph $\G$. (Recall that $\chi$ is the minimum number of independent sets ---color classes---
into which $V$ can be partitioned.)

\begin{coro}[\cite{hof70}]
Let $\G$ be a graph with  eigenvalues $\lambda_1\ge\cdots \ge\lambda_n$. Then, its chromatic
number satisfies
\begin{equation}\label{chi-vs-ind}
\chi \ge
%\frac{1}{\alpha^*} \sum_{i=1}^{\chi}\|\vecrho U_i\|^2 = \frac{\|\v\|^2}{\alpha^*}\ge
1-\frac{\lambda_1}{\lambda_n}.
\end{equation}
\end{coro}

\proof
Suppose first that
$\G$ is connected, with positive eigenvector
$\vecnu$. Since, for any minimum  coloring of $\G$, each color class $U_i$,
$1\le i\le \chi$, is an independent set, we have $\|\vecrho U_i\|^2\le
\alpha^*$. Hence, $\chi\ge  \|\vecnu\|^2/\alpha^*$ and (\ref{weight-indepen})  yields the result.
Otherwise, if $\G$ is disconnected, we only need to apply  (\ref{chi-vs-ind}) to any connected
component with maximum eigenvalue $\lambda_1$.
\endproof

A direct proof of (\ref{chi-vs-ind}) was given by Haemers \cite{ha79,ha95}.
His proof also uses eigenvalue interlacing, and so it is different from Hoffman's original one.
However, excepting for the regular case, Haemers' proof is not related to any independence-like
number. As cited by that author in \cite{ha95}, his proof  has become a common example of
application of the interlacing technique (see, for instance, Godsil
\cite[p.48]{g93} or Lov\'asz \cite[Problem 11.21]{l79-2}).

When $\G$ is regular,  Theorem \ref{theo-weight-indepen} reduces to the following bound
for the (standard) independence number:
\begin{equation}
\label{alpha}
\alpha\le \frac{n}{1-\frac{\lambda_1}{\lambda_n}}
\end{equation}
which, according to Haemers \cite{ha80,ha95}, is an unpublished result of Hoffman. The first
published proof is due to Lov\'asz \cite{l79} who derived the same upper bound
for the so-called {\it Shannon capacity} of $\G$ \cite{sha56}, defined as
$$
\Theta:=\sup_k\sqrt[k]{\alpha(\G^k)}=\lim_{k\rightarrow \infty} \sqrt[k]{\alpha(\G^k)}.
$$
Here $\alpha(\G^k)$ denotes the independence number of $\G^k$, the {\it product} of
$k$ copies of $\G$,  with vertex set
$V\times\stackrel{k}{\cdots}\times V$ and adjacencies between distinct
vertices $(u_1,\ldots,u_k)
\sim (v_1,\ldots,v_k)$ iff, for any $1\le i\le k$, either $u_i=v_i$ or $u_i\sim v_i$. Note that,
since
$\alpha(\G^k)\ge \alpha^k$, the Shannon capacity always satisfies the bound $\Theta\ge \alpha$.
For more details about this parameter, see also Knuth's paper  \cite{knu94}.
The weight version of the Shannon capacity can be defined by just writing
$$
\Theta^*:=\sup_k\sqrt[k]{\alpha^*(\G^k)}
$$
and, as expected, it can be shown to be bounded above by the weight analogue of Lov\'asz bound,
as the next theorem shows. (To prove it, recall that the {\it Kronecker product}  of two matrices
$\A\otimes \B$ is obtained by replacing each entry $(\A)_{uv}$  with the  matrix
$(\A)_{uv}\B$. Then, if $\vecnu$ and $\veceta$ are eigenvectors
of $\A$ and $\B$, with corresponding eigenvalues $\lambda$ and
$\mu$, respectively, then $\vecnu\otimes \veceta$ ---viewing $\vecnu$ and $\veceta$  as
1-column matrices--- is an eigenvector of $\A\otimes \B$, with eigenvalue
$\lambda\mu$.)

\begin{theo}
Let $\G$ be a  graph with eigenvalues $\lambda_1\ge\cdots
\ge\lambda_n$ and positive eigenvector $\vecnu$. Then, its weight Shannon capacity satisfies
\begin{equation}
\Theta^*\le \frac{\|\vecnu\|^2}{1-\frac{\lambda_1}{\lambda_n}}.
\end{equation}
\end{theo}

\proof
The proof goes along the same lines as that given by Haemers \cite{ha95} in the regular case.
As commented in Section 2,  the above results remain valid for any symmetric matrix
%%$\A^*\in \mho(\G)$, where $\mho(\G)$ is the set of (non-zero)
%%irreducible symmetric  matrices with
%%$(\A^*)_{uv}\ge 0$ if $u\sim v$, and $(\A^*)_{uv}= 0$ otherwise.
$\A^*$ with $(\A^*)_{uv}= 0$ if $u\not\sim v$, which has maximum eigenvalue with a positive
eigenvector.
%%Notice that, in this case, the
%%maximum eigenvalue $\lambda_1$ is not necessarily simple and hence the choice of an
%%associated (not necessarily positive) eigenvector $\vecnu$ is not unique
Then the application of Theorem \ref{theo-weight-indepen} to the matrix
$$
\A^*(\G^k):=(\A-\lambda_n\I)\otimes\stackrel{k}{\cdots}\otimes(\A-\lambda_n\I)-(-\lambda_n)^k,
$$
with maximum eigenvalue $(\lambda_1-\lambda_n)^k-(-\lambda_n)^k$,
positive eigenvector $\vecnu\otimes\stackrel{k}{\cdots}\otimes\vecnu$, and
minimum eigenvalue $-(-\lambda_n)^k$ gives
$$
\label{weight-k-indepen}
\alpha^*(\G^k)\le \left(\frac{\|\vecnu\|^2}{1-\frac{\lambda_1}{\lambda_n}}\right)^k,
$$
whence the result follows.
\endproof

Notice that, since  $\alpha^* \le  \Theta^*$ and $\Theta \le \Theta^*$,
%($\G$ is connected),
the above result yields also bounds for both $\alpha^*$ (that
is Theorem \ref{theo-weight-indepen}) and $\Theta$, the (standard) Shannon capacity  of a
(not necessarily regular)  graph.

%%%%%%%%%%%%%%%%%%%%%%%%%%%%%Sec. 4
\section{The weight odd-independence numbers}

The concepts of odd and even distance were introduced by Bond
and Delorme in \cite{bd88}, and they are based on looking at the  parity of the lengths
of the  walks considered. Thus, the {\it odd distance} between two (not necessarily
different) vertices $u,v$ of a graph $\G$, denoted by $\partial_\o(u,v)$, is the length of a
shortest walk of odd length between them. By using odd distances, we can now consider other
related metric parameters, such as the {\it odd diameter} $D_\o$ and the {\it odd girth} $g_\o$,
defined as expected (in $D_\o$ we must also consider $\partial_\o(u,u)$ when looking at
maximum odd distance between pairs of vertices). Thus, since
$\G$ has no loops,
$1<g_\o\le
\d_\o(u,u)\le D_\o$ for any vertex $u\in V$ and, if $\G$  bipartite,
$g_\o=D_\o=\infty$. Otherwise, the above-mentioned authors proved that
$D_\o\le 2D+1$, with $D$ being the standard diameter of $\G$. In fact it can be
shown that, if $\G$ is a non-bipartite connected graph, then $D_\o\le d^\star\le 2d+1$, where
$d^\star$ is the number of points of the ``symmetrized mesh"
$\M^\star:=\M\cup\{0\}\cup(-\M)$, with $\M=\ev \G \setminus\{\lambda_1\}$, and $d=|\M |$ (see
\cite{f1}).

By using odd distances, the author \cite{f1} introduced a new measure of independence
as follows. Let $k\ge 1$ be an odd integer. Then,  the {\it odd-$k$-independence number}
$\alpha_{\o k}$ is defined as the maximum number of vertices which are at odd distance greater
than $k$ from each other (including the odd distance from one vertex to itself). Note that, in
particular, $\alpha_{\o 1}$ coincides with the independence number $\alpha$. Any set of vertices
satisfying such a condition is called an {\it odd-$k$-independent set}, so that the corresponding
weight parameter is
$$
\alpha_{\o k}^*:=\max_{U\subset V}\{\|\vecrho U\|^2: U \mbox{ is an odd-$k$-independent set}\}.
$$
Basically the same proof used in Theorem \ref{theo-weight-indepen} yields the following
result, whose first part was also proved in \cite{f1} by using another technique.

\begin{theo}
\label{theo-odd-indepen}
Let $\G$ be a graph with eigenvalues $\lambda_1\ge \lambda_2\ge\cdots\ge\lambda_n$, and
positive eigenvector $\vecnu$. Let $q$ be a polynomial with only odd powers, degree $k$,
$q(\lambda_1)>0$,  and $q_{\min}:=\min_{2\le i\le n}\{q(\lambda_i)\}$. Then, provided that
$\alpha_{\o k}^*\neq 0$, we have
\begin{equation}
\alpha_{\o k}^* \le \frac{\|\vecnu\|^2}{1-\frac{q(\lambda_1)}{q_{\min}}}.
\end{equation}
If the bound is attained for some odd-$k$-independent set $C$, then
\begin{equation}
q(\A)\vecrho C =-q_{\min}\vecrho \overline{C}.
\end{equation}
\end{theo}

\proof
Let $C\subset V$ be an odd-$k$-independent set with
$\alpha_{\o k}^*=\|\vecrho C\|^2$. By the hypotheses on the polynomial $q$, the matrix $q(\A)$
has minimum eigenvalue $\zeta_n:=\min\{q(\lambda_1),q_{\min}\}$ and  $(q(\A))_{uv}=0$ for any
$u,v\in C$.   Hence,  the weight-quotient matrix of $q(\A)$, denoted by $\B_q$, with respect to
the partition $V=C\cup \Cb$ is
\begin{equation}
\label{B_q}
\B_q=\left(\begin{array}{cc}
0 & \frac{\|\vecrho C\|} {\|\vecrho\Cb\|} \\
\frac{\|\vecrho C\|} {\|\vecrho\Cb\|} &
\frac{\|\vecrho\Cb\|^2-\|\vecrho C\|^2}{\|\vecrho\Cb\|^2}
\end{array}\right)q(\lambda_1)
\end{equation}
(compare with the matrix $\B$ in (\ref{B})), with eigenvalues $\mu_1=q(\lambda_1)$  and
$\mu_2 =\frac{-q(\lambda_1)\|\vecrho C\|^2}{\|\vecnu\|^2-\|\vecrho C\|^2}$ satisfying
\begin{equation}
\label{basic-inequality}
0>\frac{-q(\lambda_1)\|\vecrho C\|^2}{\|\vecnu\|^2-\|\vecrho C\|^2}\ge \zeta_n=q_{\min},
\end{equation}
where we have used again Lemma \ref{interlac-coro}, and the hypotheses $q(\lambda_1)>0$,
$\|\vecnu\|>\|\vecrho C\|>0$. Hence the first statement follows.
Furthermore, if we get equality for some set $C$, we
have $q(\lambda_1)=-q_{\min} \frac{\|\vecrho\Cb\|^2}{\|\vecrho C\|^2}$ and the interlacing is
tight: $\S\B_q=q(\A)\S$. But $\S$ consists of the two (column)  vectors
$\frac{1}{\|\vecrho C\|}\vecrho C$ and $\frac{1}{\|\vecrho \Cb\|}\vecrho \Cb$, so that the above
matrix equation reads:
$$
\left(\begin{array}{cc}
\frac{1}{\|\vecrho C\|}\vecrho C  & \frac{1}{\|\vecrho \Cb\|}\vecrho\Cb
\end{array}\right)
\left(\begin{array}{cc}
0 & \frac{\|\vecrho\Cb\|}{\|\vecrho C\|} \\
\frac{\|\vecrho\Cb\|}{\|\vecrho C\|} & \frac{\|\vecrho\Cb\|^2-\|\vecrho C\|^2}{\|\vecrho C\|^2}
\end{array}\right)(-q_{\min})
=q(\A)\left(\begin{array}{cc}
\frac{1}{\|\vecrho C\|}\vecrho C  & \frac{1}{\|\vecrho \Cb\|}\vecrho\Cb
\end{array}\right)
$$
giving $\vecrho \Cb (-q_{\min})=q(\A)\vecrho C$, as claimed.
\endproof

From the above proof, note that, if $q_{\min}\ge 0$, then
(\ref{basic-inequality}) gives a contradiction and hence it must be $\alpha_{\o k}^*= 0$. This
implies the existence of an odd closed walk of length at most $k$ through any vertex and,
therefore, $g_\o\le k$. From these facts, we easily deduce that, if $\G$ is not bipartite
($\lambda_n\neq-\lambda_1$), then $g_\o\le d^\star$, where $d^\star$ is the number of points
of the symmetrized mesh $\M^\star$ defined above. (Just let
$q$ be any polynomial having such points as its roots and taking positive value at $\lambda_1$.)

As another consequence of Theorem \ref{theo-odd-indepen}, and reasoning as in the previous
section, we can now derive an upper bound for a chromatic-like number, which we could call the
``odd-$k$-chromatic number." Let $\G$ be a graph with $n$ vertices and odd girth $g_\o$. For
each odd integer $k$, $1\le k< g_\o$,  the {\it odd-$k$-chromatic number} of $\G$, denoted by
$\chi_{\o k}=\chi_{\o k} (\G)$, is the minimum number of colors  that can be assigned to the
vertices of $\G$ in such a way that any two (not necessarily different) vertices having the same
color are at odd distance greater than $k$ from each other. Notice that, with this definition,
$$
\chi_{\o 1}(\equiv \chi)\le
\chi_{\o 3} \le\cdots \le \chi_{\o g_\o-2}\le n
$$
(if $\G$ is bipartite, $\chi_{\o k}=2$ for any
odd $k\ge 1$).  In other words, we can say that $\chi_{\o k}$ is the
minimum number of odd-$k$-independent sets into which $V$ can be partitioned. Within this
framework, the following result could be seen as a generalization of Hoffman's bound
(\ref{chi-vs-ind}).

\begin{coro}
\label{odd-chromatic}
Let $\G$ be a graph with  odd girth $g_\o$ and eigenvalues
$\lambda_1\ge\lambda_2\ge\cdots\ge\lambda_n$. Let $q$ be
a polynomial of degree $k$ as above. Then, for any odd integer $k$, $1\le k< g_\o$,  the
odd-$k$-chromatic number satisfies
\begin{equation}
\label{k-chromatic}
\chi_{\o k} \ge 1-\frac{q(\lambda_1)}{q_{\min}}.
\end{equation}
\end{coro}

\proof
Since $k<g_\o$, we have $\alpha^*_{\o k}\ge 1$. Then Theorem \ref{theo-odd-indepen}
applies and the result follows from $\chi_{\o k}\ge \|\vecnu\|^2/\alpha^*_{\o k}$.
\endproof

In particular, if we take $q(x)=(\frac{x}{-\lambda_n})^k$, with $q_{\min}=-1$, we get
\begin{equation}
\label{chi-k}
\chi_{\o k} \ge 1-\left(\frac{\lambda_1}{\lambda_n}\right)^k,
\end{equation}
a result to be compared with  (\ref{chi-vs-ind}).

Of course, we can do better if we look for the (odd) polynomials, with  degree at most $k$, that
maximize the quotient
$q(\lambda_1)/(-q_{\min})$ or, alternatively, we can try to maximize $q(\lambda_1)$
among the polynomials  that satisfy $q_{\min}\ge -1$. These polynomials were studied with
some detail in \cite{f1}. Also, a method to compute them, based on solving a linear programming
problem, was proposed. They will be referred to as the odd polynomials and
denoted by $Q_k$. To discuss some of their properties, it is better to consider only the
distinct eigenvalues of the graph: $\ev \G=\{\theta_0>\theta_1>\cdots>\theta_d\}$. As before,
set $\M =\ev \G \setminus\{\theta_0\}$ and consider  the symmetrized mesh
$\M^\star=\M\cup\{0\}\cup(-\M)$, with  $d^{\star}:=|\M^{\star}|$ points.
Then, for any  odd integer $k$, $1\le k \le d^{\star}-2$, the {\it odd polynomial} $Q_k$ satisfies
\begin{equation}
\label{def-Q_k}
Q_k(\theta_0) = \max_{q\in \R^-_{k} [\theta]} \{ q(\theta_0) : q_{\min} \ge -1 \}
\end{equation}
where $\R^-_{k} [\theta]$ stands  for the set of real
polynomials with only odd powers and degree at most $k$, and $q_{\min}=\min_{\theta\in \M
}\{q(\theta)\}$. If $\G$ is not  bipartite ($\pm\theta_0\not\in \M^\star$) it was shown that
there is a unique odd polynomial of degree $k$,  satisfying ${(Q_k)}_{\min}=-1$, and
\begin{equation}
\label{increasing-odd-pol}
Q_1(\theta_0)< Q_3(\theta_0)< \cdots <Q_{d^\star-2}(\theta_0).
\end{equation}
%(If $\G$ is bipartite,  $Q_k(\theta_0)=1$ for any $k$.)
In particular, the extremal cases $k=1$ and $k=d^{\star}-2$ admit
closed expressions. Namely,
$Q_1(x)=\frac{x}{-\theta_d}$, and
$Q_{d^{\star}-2}$ being the polynomial which takes alternating values $\pm 1$ at
\linebreak
$\M^\star\setminus\{ 0\}=\{\vartheta_1>\vartheta_2>\cdots>\vartheta_{d^\star} \}$. This
gives (using Lagrange interpolation in the second case):
\begin{equation}
\label{Q(theta)}
Q_1(\theta_0)=\frac{\theta_0}{-\theta_d}, \ \ \
Q_{d^{\star}-2}(\theta_0)=\sum_{i=1}^{d^\star}\frac{\pi_0}{\pi_i}
\end{equation}
where $\pi_i:=\prod_{j=0,j\neq i}|\vartheta_i-\vartheta_j|$\ \  ($\vartheta_0=\theta_0$).

Then, in terms of these polynomials and using the new notation for the eigenvalues, Theorem
\ref{theo-odd-indepen} reads
\begin{equation}
\label{opt-odd-indepen}
\alpha_{\o k}^* \le \frac{\|\vecnu\|^2}{1+Q_k(\theta_0)}
\end{equation}
and, in the case of equality for some vertex subset $C$,
\begin{equation}
Q_k(\A)\vecrho C = \vecrho\overline{C};
\end{equation}
whereas Corollary
\ref{odd-chromatic} yields
\begin{equation}
\label{Q-k-chromatic}
\chi_{\o k} \ge 1+Q_k(\theta_0).
\end{equation}

\begin{example}
Let $\G=O_4$, the  (regular) {\it ``odd graph" } with degree $4$,  $35$ vertices, and eigenvalues
$\ev O_4 =\{4>2>-1>-3\}$. (The odd graph $O_k$ has  the $(k-1)$-subsets of a
$(2k-1)$-subset as vertices, and two vertices are adjacent iff their corresponding subsets are
disjoint; see Biggs
\cite{b79,b93}.) Then the corresponding symmetrized mesh is $\M^{\star}=\{ 0,\pm 1, \pm 2,
\pm 3\}$ and hence
$g_\o\le D_\o\le 7$ (in fact $g_\o= 7$). The corresponding odd polynomials and their values at
$\theta_0=4$ are:
\begin{itemize}
\item
$Q_5(x)=\frac{1}{12}\left(x^5-11 x^3+22x \right)$,\  \ $34$;
\item
$Q_3(x)=\frac{1}{6}\left( x^3-7x\right)$,\ \   $6$;
 \item
$Q_1(x)=\frac{1}{3} x$,\ \  $4/3$.
\end{itemize}
\noindent
Hence, the respective bounds for the odd-$k$-independence numbers, given by
$(\ref{opt-odd-indepen})$, turn out to be $\alpha_{\o 1}=\alpha\le 15$, $\alpha_{\o 3}\le 5$, and
$\alpha_{\o 5}\le  1$.  In fact, all these bounds are tight, as can be easily shown by  using the
known formulas for the distances between vertices in the odd graphs (see Biggs \cite{b79}).
\end{example}

When $\G$ is bipartite we have $Q_k(\theta_0)=-Q_k(\theta_d)=1$ for any $k$, and
(\ref{opt-odd-indepen}) yields  $\alpha_{\o k}^*\le \|\vecnu\|^2/2$, as expected. In
the case of regular non-bipartite connected graphs, and since
$\alpha_{\o g_\o-2}\ge 1$, (\ref{opt-odd-indepen}) gives the following result.

\begin{coro}
The order $n$ of a non-bipartite regular connected graph $\G$, with eigenvalues
$\ev \G=\{\theta_0>\theta_1>\cdots>\theta_d\}$ and odd girth $g_\o$, satisfies the bound
\begin{equation}
\label{order-bound}
n\ge  Q_{g_\o-2}(\theta_0) + 1
\end{equation}
where $Q_{g_\o-2}$ is the odd $(g_\o-2)$-polynomial.
\end{coro}
From the example above, note that the bound (\ref{order-bound}) is tight for $O_4$. In fact,
using the value of $Q_{d^\star-2}(\theta_0)$ given in (\ref{Q(theta)}), it can be shown that this
property is shared by all odd graphs $O_k$. Notice also that (\ref{order-bound}) still holds if
we replace $g_\o$ by the standard girth $g$ (since $g_\o\ge g$ and the odd polynomials satisfy
(\ref{increasing-odd-pol})).
%%$Q_1(\theta_0)<Q_3(\theta_0)<\cdots<Q_{d^\star-2}(\theta_0)$).

The next straightforward consequence of Theorem \ref{theo-odd-indepen} is also given in terms
of the odd polynomials. Let $\G$ be a graph on $n$ vertices. Given any integer $1\le t\le n$,
let us define the {\it odd $t$-diameter} of a graph $\G$ as
$$
D_{\o t}:=\max_{U\subset V} \{\min_{u,v\in U}\partial_\o (u,v):  |U|=t\}
$$
so that the following inequalities hold
$$
D_\o\ge D_{\o 1}\ge D_{\o 2}\ge \cdots\ge D_{\o n}(=1)
$$
where $D_\o$ is the above-mentioned odd diameter.

\begin{coro}
Let $\G$ be a graph as above, with odd polynomials $Q_k$. Then,
\begin{equation}
 Q_k(\theta_0)> \frac{\|\vecnu\|^2}{t}- 1\ \  \Rightarrow\ \   D_{\o t} \le k.
\end{equation}
\end{coro}

\proof
Under the hypothesis, (\ref{opt-odd-indepen}) gives $\alpha_{\o k}^*<t$. Consequently,
between any $t$ vertices, there must be some walk of odd length  $\le k$ (perhaps between a
vertex and itself).
\endproof

%%%%%%%%%%%%%%%%%%%%%%%%%%%%Sec. 5
\section{The weight set independence numbers}

In this section we study another generalization of the concept of independence, which concerns
the elements considered (sets instead of single vertices) rather than the type of distance
involved.  Indeed, we can extend the notion of $k$-independence to  vertex subsets if
we require that they must be at distance greater than $k$ from each other. We will first suppose
that all such subsets have the same weight. Afterwards, we shall pay attention to the simplest
case of (two) subsets with different weights.

\subsection{Subsets with equal weights}
Assume that the graph $\G$, on $n$ vertices, has some vertex subset $U$ with weight
$w:=\|\vecrho U\|^2$. (Note that, if $U=\{u\}$, the ``weight" of vertex $u$ is now $\nu_u^2$.)
Then, given some integer
$k\ge 0$,
%%and real number $1\le w\le \|\vecnu\|^2$,
we define the {\it $(w,k)$-independence number}, denoted by
$\alpha_k^w$, as the maximum number of $k$-independent subsets
%%$, if any,
with common weight $w$.
%%If there is no subset with weight $w$ it will be convenient to define
%%$\alpha_k^w:=0$. Otherwise,  and
As in the standard notion of independence, we assume that the set
$U$ is $k$-independent from itself, so that $\alpha_k^w\ge1$.
Notice that, when $\G$ is regular, $\alpha_k^1$ is the maximum number of vertices which are
mutually at distance greater than $k$. This parameter, denoted just  by $\alpha_k$, has been
recently considered in the literature by Delorme and  Tillich \cite{dt97}, and Garriga, Yebra
and the author \cite{f3,fg1,g97}, and it is called the  {\it
$k$-independence number}. Thus, $\alpha_0=n$, $\alpha_1\equiv \alpha$, and $\alpha_k$ is, in
fact, the independence number of the $k$-th power of $\G$ (that is the graph $\G_{\le k}$ with
vertex set
$V$ and where two vertices are adjacent whenever their distance in $\G$ is at most $k$).
In order to give bounds for $\alpha_k$, it is useful to consider the so-called alternating
polynomials, introduced in \cite{fgy1}, which can be thought of as the discrete version of the
Chebychev polynomials.
As above, let $\G$ be a graph with $\ev\G=\{\theta_0>\theta_1>\cdots>\theta_d\}$. For
any integer $0\le k\le d-1$, the {\it $k$-alternating polynomial}
$P_k$,  is  the (unique) polynomial of degree $k$ satisfying
\begin{equation}
\label{def-P_k}
P_k (\theta_0)=\max_{p\in \R_k[x]}\{p(\theta_0):\, \|p\|_\infty\le 1\}
\end{equation}
where $\|p\|_\infty:=\max_{1\le i\le d}|p(\theta_i)|$.
%%When $k=d-1$, we simply speak about the {\it alternating polynomial} and denote it by
%%$P:=P_{d-1}$.
Thus, we obviously have  $P_0=1$. Otherwise, for $k\ge 1$, it was proved in
\cite{fgy1} that the
$k$-alternating polynomial is characterized by taking $k+1$
alternating values $\pm 1$ at  $\ev \G\setminus \{\theta_0\}$, with
$P_k(\theta_1)=1$ and $P_k(\theta_d)=(-1)^k$. Moreover,
\begin{equation}
\label{increasing-alt-pol}
1<P_1(\theta_0)< P_2(\theta_0)< \cdots <P_{d-1}(\theta_0).
\end{equation}
In particular, for the values
$k=1$ and $k=d-1$, the above characterization gives
$P_1(x)=2\frac{x-\theta_1}{\theta_1-\theta_d}+1$, and the $(d-1)$-alternating polynomial is
defined by $P_{d-1}(\theta_{i})=(-1)^{i+1}$,
$1\le i\le d$. Thus, (using again Lagrange interpolation) we get
\begin{equation}
\label{P(theta)}
P_1(\theta_0)=2\frac{\theta_0-\theta_1}{\theta_1-\theta_d}+1, \ \ \
P_{d-1} (\theta_0)=\sum_{i=1}^d \frac{\pi_0}{\pi_i}
\end{equation}
where $\pi_i:=\prod_{j=0,j\neq i}^d|\theta_i-\theta_j|$, $0\le i\le d$.
Some particular cases of these polynomials were also considered by Van Dam and Haemers in
\cite{vdh95}. In fact, as noted by Van Dam \cite{vd96}, they had already been considered in the
theory of uniform approximations of continuous functions.
%%%%%%%%%%%%%%%%%%%%%

In terms of the alternating polynomials, the author  \cite{f3} showed
 that, for a regular connected graph $\G$ on $n$ vertices,
 %%and  $\ev \G=\{\theta_0> \theta_1> \cdots> \theta_d\}$
the $k$-independence number is  bounded above by
\begin{equation}
\label{k-indepen}
\alpha_k \le \frac{2n}{P_k (\theta_0)+1}.
\end{equation}
In the next theorem the above result is generalized by giving a similar bound
for $\alpha_k^w$ of any (connected) graph. The case $w>\|\vecnu\|^2/2$  can be
excluded  since then $\alpha_k^w=1$ for any $k$.

\begin{theo}
\label{theo-(w,k)-indepen}
Let $\G$ be a connected graph with eigenvalues $\ev \G=\{\theta_0> \theta_1> \cdots
> \theta_d\}$, positive eigenvector $\vecnu$, and $k$-alternating polynomials $P_k$, $0\le k\le
d-1$.  Assume that, for some weight $w\ge 1$, the $(w,k)$-independence number satisfies
$\alpha_k^w\ge 2$. Then,
\begin{equation}
\label{(w,k)-indepen}
\alpha_k^w \le  \frac{2\|\vecnu\|^2}{w(P_k(\theta_0)+1)}.
\end{equation}
\end{theo}

\proof
We can suppose that $k\ge 1$ since, otherwise, $P_0=1$ and the result trivially
holds. Then, let $r:=\alpha_k^w<\|\vecnu\|^2/w$, and assume that $U_i$,
$1\le i\le r$, are some  $k$-independent sets  with common weight $w=\|\vecrho U_i\|^2$.  Take
the polynomial $q:=\frac{r}{2}P_k+\frac{r-2}{2}$, which satisfies $-1\le q(\theta_i) \le r-1$ for
any $1\le i\le d$. Then, as $P_k(\theta_0)> 1$ and $\G$ is connected, the  matrix $q(\A(\G))$ has
eigenvalues $\{-1<\cdots< r-1< q(\theta_0)\}$ and $q(\theta_0)$ has
multiplicity $1$. Moreover,  the complete graph $K_r$ with vertex set
$\{1,2,\ldots,r\}$ has eigenvalues $\ev K_r=\{-1< r-1\}$. Consequently, the matrix
$\K$ obtained as the Kronecker product $\A(K_r)\otimes q(\A(\G))$ has eigenvalues
$$
\{-q(\theta_0)< -(r-1)< \cdots <(r-1)^2< (r-1)q(\theta_0)\}.
$$
Let us now consider the partition $\bigcup_{i=1}^r [(i,U_i)\cup(i,\overline{U_i})]$ of
(the rows and columns of) $\K$. Since, for $i\neq j$,
$\partial(U_i,U_j)>k$ in $\G$, the weight-quotient matrix with respect to such a partition
turns out to be again a Kronecker product, namely $\B:=\A(K_r)\otimes \B_q$, with
$\B_q$  as in (\ref{B_q}) ---that is,  the weight-quotient matrix of $q(\A(\G))$ with respect to
any partition $U_i\cup \overline{U_i}$---  with eigenvalues $\{-q(\theta_0)w/(\|\vecnu\|^2-w)<
q(\theta_0)\}$. Therefore, $\B$ has eigenvalues
$$
\textstyle
\{-q(\theta_0)< -\frac{(r-1)w}{\|\vecnu\|^2-w}q(\theta_0)<
\frac{w}{\|\vecnu\|^2-w}q(\theta_0)< (r-1)q(\theta_0)\}.
$$
Thus, since for both matrices $\K$ and  $\B$ the minimum eigenvalue
$-q(\theta_0)$ has multiplicity $r-1$ (that is the multiplicity of $-1$ as eigenvalue of
$\A(K_r)$), we have, by Lemma \ref{interlac-coro}, that their $r$-th smallest eigenvalues
satisfy
$$
-(r-1)\le -\frac{(r-1)w}{\|\vecnu\|^2-w}q(\theta_0).
$$
Hence, using the expression for $q$,
\begin{equation}
\label{(w,k)-indepen-bis}
 \frac{r}{2}\left(P_k(\theta_0)+1\right)-1=q(\theta_0)\le \frac{\|\vecnu\|^2}{w}-1
\end{equation}
 whence we get (\ref{(w,k)-indepen}).
%and solving for $r$ the theorem follows.
\endproof

As in the previous section, let us now give some straightforward consequences of the above
theorem.  First,  from its  proof we get the following simple corollary.

\begin{coro}
\label{coro-(w,k)-indepen}
Assume that the connected graph $\G$ has a vertex subset $U$ with weight $w_U\ge w$,
for some $1\le w <\|\vecnu\|^2$. Then, if
$
P_k(\theta_0) >\frac{\|\vecnu\|^2}{w}-1
$,
all the other subsets with weight $w_U$, if any, are at distance at most $k$ from $U$.
\end{coro}

\proof
From the hypotheses, we  have $P_k(\theta_0)
>\frac{\|\vecnu\|^2}{w_U}-1$. Consequently, (\ref{(w,k)-indepen-bis}) gives
$r< 2$, a contradiction. Hence, it must be $\alpha_k^{w_U}= 1$ and the result follows.
\endproof
%%In particular, since $w_u\ge 1$ for any vertex $u\in V$, we obtain the following  result
%%given in \cite{fgy1}:
%%$$
%%P_k(\theta_0) >\|\vecnu\|^2-1 \ \ \Rightarrow \ \ D(\G) \le k.
%%$$

%%For simplicity, we first  limit ourselves to the regular case.
Consider now the  specialization of the above results to regular graphs.
In this case,  for any integer $w$,
$1\le w<n$, we can consider subsets of any weight (cardinality) $w$. Then,  with the notation
$\alpha_k\equiv \alpha_k^1$, we clearly have $\alpha_k^w\ge \lfloor\alpha_k/w\rfloor$.
Moreover, Theorem \ref{theo-(w,k)-indepen}
%%and Corollary \ref{coro-(w,k)-indepen}
gives:

\begin{coro}
\label{coro-k-indepen}
Let $\G$ be a $\delta$-regular connected graph on $n$ vertices,
$\ev \G=$ $\{\theta_0(=\delta)>\theta_1>\cdots >\theta_d\}$, and
$k$-alternating polynomials $P_k$, $0\le k\le d-1$.  Then, for any integer $w$, $1\le w< n$, the
$(w,k)$-independence number satisfies
\begin{equation}
\label{mk-indepen}
\alpha_k^w \le \max\left\{1,\frac{2 n}{w(P_k (\theta_0)+1)}\right\}.
\end{equation}
\end{coro}

In particular,  if $\G$ contains at least two
($1$-)independent $w$-sets, then  $\alpha_1^w$ satisfies the bound
\begin{equation}
\label{alpha^m}
\alpha_1^w\le \frac{n(\theta_1-\theta_d)}{w(\theta_0-\theta_d)}
\end{equation}
where we have used the value of $P_1(\theta_0)$ in (\ref{P(theta)}).
Notice that, for all non-complete connected graphs with $\theta_1>0$ (that is, those different
from the complete multipartite graphs), the bound for the independence number $\alpha(\ge 2)$
obtained by taking  $w=1$ in (\ref{alpha^m}) is worse than $\alpha \le
\frac{n(-\theta_d)}{\theta_0-\theta_d}$, given in (\ref{alpha}).
Another  particular case of (\ref{mk-indepen}) worth mentioning is the following upper bound
for the maximum number of $w$-sets of a regular connected graph which are pairwise at
(spectrally maximum) distance $d$.
$$
\alpha_{d-1}^w\le \max\left\{1,\frac{2n}{w\sum_{i=0}^{d}\frac{\pi_0}{\pi_i}}\right\}
$$
where we have used the value of $P_{d-1}(\theta_0)$ in (\ref{P(theta)}).
%$\pi_i:=\prod_{j=0, j\neq i}^d|\theta_i-\theta_j|$.

%%\begin{coro}
%%Let $\G$ be a regular graph as above.
%%containing some vertex set $U$ with $w=\|\vecrho U\|^2$.
%%Then, its $(m,k)$-chromatic number satisfies
%%\begin{equation}
%%\chi_k^m\ge \frac{m}{2}\left( P_k(\theta_0)+1\right).
%%\end{equation}
%%\end{coro}

Assume now that the (not necessarily regular) graph $\G$ has at least $t\ge 2$ vertex subsets
$U_1,\ldots,U_t$ with the same weight $w$, say. Then,  we  can define the {\it $(w,t)$-diameter}
by
\begin{equation}
\label{wk-diameter}
D_t^w:= \max_{U_1,\ldots,U_t\subset V}\{  \min_{1\le i<j\le t}
\partial(U_i,U_j) : \|\vecrho U_i\|^2=w, 1\le i\le t \}.
\end{equation}

Thus, if $\G$ is regular, the $(w,t)$-diameter coincides with the parameter $D_{t\times w}$,
studied by  Garriga and the author in \cite{fg1} (there we consider the minimum distance
between families of $t$ subsets on $w$ vertices). Similarly, the $t$-diameter $D_t$ considered
by Chung,  Delorme, and  Sol\'e \cite{cds97},  corresponds to
$D_t^1$ (that is, the largest pairwise minimum distance of a set of $t$ vertices), so that the
(standard) diameter is just $D_2^1$.
% (Eigenvalue upper bounds for $D_k^1$ were given in...)

Now,
%reasoning as in Corollary \ref{},
Theorem \ref{theo-(w,k)-indepen} gives the following
result which can be seen as an extension of Corollary \ref{coro-(w,k)-indepen} (the case $t=2$).

\begin{coro}
Let $\G$ be a graph as above containing at least $t\ge 2$ vertex subsets $U$ with weight
$w=\|\vecrho U\|^2$. Then,
\begin{equation}
\label{coro-wt-diam}
P_k (\theta_0) > \frac{2\|\vecnu\|^2}{wt}-1 \ \  \Rightarrow\ \  \ D_t^w\le k .
\end{equation}
\end{coro}

When $\G$ is a regular graph on $n$ vertices (with ${n\choose w}\ge t$) we  have the following
result concerning  subsets of $w$ vertices:
\begin{equation}
\label{reg-wt-diam}
P_k (\theta_0) > \frac{2n}{wt}-1 \ \  \Rightarrow \ \  \ D_t^w\le k.
\end{equation}
The particular case  $w=1$ was  proved in \cite{f3} by using a different technique.
Notice that the above results still hold if we replace $P_k$ by the Chebychev polynomial
$T_k$  ``shifted" from
$[-1,1]$ to $[\theta_d,\theta_1]$, that is
$T_k^*(x):=T_k\left( \frac{2x-\theta_1-\theta_d}{\theta_1-\theta_d} \right)$ (since
$\|T_k^*\|_\infty=\|P_k\|_\infty=1$ and $P_k(\theta_0)\ge T_k^*(\theta_0)$). Then, using that
$T_k(x)=\cosh (k\cosh^{-1} x)$, (\ref{reg-wt-diam}) yields:
\begin{equation}
D_t^w \leq \left\lfloor\frac{{\rm cosh}^{-1} \left(\frac{2n}{wt}-1\right)}
 {{\rm cosh}^{-1}\left(\frac{2\theta_0-\theta_1
- \theta_d}{\theta_1-\theta_d} \right)} \right\rfloor +1.
\end{equation}
A result to be compared with that given by Kahale  \cite{k97}, who proved that if $\G$ is a
regular connected graph on $n$ vertices, and
$\vartheta_1(=\theta_0),\vartheta_2,\ldots,\vartheta_n$
represent its eigenvalues with absolute value in non-increasing order,
$|\vartheta_1|>|\vartheta_2|\ge\cdots\ge |\vartheta_n|$, then
\begin{equation}
D_t^w \leq \left\lceil\frac{{\rm cosh}^{-1} \left(\frac{n}{w}-1\right)}
{{\rm cosh}^{-1} ({\vartheta_1 / |\vartheta_t|})} \right\rceil +1.
\end{equation}

\subsection{Subgraphs with different weights}

When we consider vertex subsets with different weights, we can still apply the same
techniques as above. However the complexity of the analysis
steadily (dramatically) increases with the number of sets considered. By way of example, we
analyze below the simplest case of two subsets. In this context, the following theorem
was also proved in \cite{fg3}  without  using eigenvalue interlacing. The corresponding
results for either regular graphs or Laplacian spectrum were also proved by Van Dam and
Haemers \cite{vdh95} and Van Dam \cite{vd98}, respectively.

\begin{theo}
Let $\G$ be a connected graph with eigenvalues $\ev \G=\{\theta_0>\theta_1>\cdots>\theta_d\}$
and
$k$-alternating polynomials $P_k$, $0\le k\le d-1$. Let
$X,Y$ be  two subsets of vertices such that
$\partial(X,Y)>k$. Then,
$$
P_k(\theta_0)\le \frac{\|\vecrho \Xb\|\|\vecrho \Yb\|}
{\|\vecrho X\|\|\vecrho Y\|}.
$$
\end{theo}

\proof
Since $\partial(X,Y)\ge 1$, we have $\|\vecrho X\|\le
\|\vecrho \Yb\|$ and $\|\vecrho Y\|\le \|\vecrho \Xb\|$, and hence the result is trivial for $k=0$.
The proof for $k\ge 1$ is similar to that of Theorem  \ref{theo-(w,k)-indepen}, but taking
$r=2$. Then, the polynomial $q$ is just $P_k$ and hence we  consider the matrix
$\K :=\A(K_2)\otimes P_k(\A(\G))$ with eigenvalues
$\pm P_k(\theta_0),\ldots, \pm P_k(\theta_d)$, satisfying $|P_k(\theta_i)|\le
\|P_k\|_\infty=1< P_k(\theta_0)$, for $1\le i\le d$, and  $P_k(\theta_0)$ having multiplicity
$1$.  Moreover, the weight-quotient matrix of $\K$, with respect to the  partition
$(1,X)\cup(1,\Xb)\cup(2,Y)\cup(2,\Yb)$, is now
$$
\B  =\left(
\begin{array}{cccc}
0 & 0 & 0 & \frac{\|\vecrho X\|}{\|\vecrho \Yb\|} \\
0 & 0 & \frac{\|\vecrho Y\|}{\|\vecrho \Xb\|}
& \frac{\|\vecrho \Xb\|^2-\|\vecrho Y\|^2}{\|\vecrho \Xb\|\|\vecrho \Yb\|} \\
0 & \frac{\|\vecrho Y\|}{\|\vecrho \Xb\|} & 0 & 0 \\
\frac{\|\vecrho X\|}{\|\vecrho \Yb\|} & \frac{\|\vecrho \Yb\|^2-\|\vecrho X\|^2}{\|\vecrho
\Yb\|\|\vecrho
\Xb\|} & 0 & 0
\end{array}
\right) P_k(\theta_0)
$$
with eigenvalues $\pm P_k(\theta_0)$ and $\pm P_k(\theta_0)\frac{\|\vecrho X\|\|\vecrho
Y\|}{\|\vecrho \Xb\|\|\vecrho \Yb\|}$.  Hence, the result follows from
%\linebreak
$$
P_k(\theta_0)\frac{\|\vecrho X\|\|\vecrho Y\|}{\|\vecrho
\Xb\|\|\vecrho \Yb\|}\le \|P_k\|_\infty=1.\quad \Box
$$

Some consequences of this theorem, together with the study of the case in which equality is
attained, can be found in \cite{fg3,fgy5}. For instance, a straightforward reasoning gives an
upper bound for the so-called {\it conditional  $(s,t)$-diameter} of $\G$, defined in
\cite{bcff95} by
$$
D_{(s,t)}=\max_{U_1,U_2\subset V}\{\partial(U_1,U_2) : |U_1|=s,\,|U_2|=t\}
$$
for some integers $1\le s,t\le n$. Namely,
$$
P_k(\theta_0) > \sqrt{\left(\frac{\|\vecnu\|^2}{s}-1\right)\left(\frac{\|\vecnu\|^2}{t}-1\right)}
\ \ \Rightarrow \ \ D_{(s,t)}\le k.
$$
(See also \cite{fgy1} for the case $s=t=1$ corresponding to the standard  diameter.)
Using weights instead of cardinalities we also get another generalization of Corollary
\ref{coro-(w,k)-indepen} which, roughly speaking, tell us that vertex subsets with large weight
tend to be close together. More precisely,  if $P_k(\theta_0) >\|\vecnu\|^2/w-1$, then all
subsets with weight at least $w$ are at most $k$ apart from each other.

%%%%%%%%%%Sec. 6
\section{Subgraphs}

For a given integer $k\ge 0$, let us consider the graph $\G_{>k}$ with the same vertex set as
$\G$, and where two vertices are adjacent iff their distance apart in $\G$ is greater than $k$. In
other words, if $k\ge1$, $\G_{>k}$ is the complement of $\G_{\le k}$, the $k$-th power of $\G$,
and $\G_{>0}=K_n$. The next result can be seen as a generalization of Theorem
\ref{theo-(w,k)-indepen}.

\begin{theo}
\label{theo-ineq-subgraphs}
Let $\G$ be a connected graph with  eigenvalues  $\ev
\G=\{\theta_0> \theta_1> \cdots > \theta_d\}$, positive eigenvector $\vecnu$, and
$k$-alternating polynomials
$P_k$, $0\le k\le d-1$. Let $H$ be a (non-trivial) subgraph of the complement of the $k$-th
power of
$\G$, $H\subseteq \G_{>k}$, with all its vertices $u$ having equal weight $w=\nu_u^2\le
\|\vecnu\|^2/2$ (in $\G$), and eigenvalues $\ev H=\{\eta_0> \eta_1> \cdots> \eta_e\}$. Then,
\end{theo}
\begin{equation}
\label{ineq-subgraphs}
1-\frac{\eta_0}{\eta_e}\le \frac{2\|\vecnu\|^2}{w(P_k(\theta_0)+1)}.
\end{equation}

\proof
Reason as in the proof of Theorem \ref{theo-(w,k)-indepen}, but
using now the polynomial
$q:=\frac{\eta_0-\eta_e}{2}P_k+\frac{\eta_0+\eta_e}{2}$, with $\eta_e\le q(\theta_i) \le
\eta_0\le q(\theta_0)$, $1\le i\le d$. Moreover, suppose that $k\ge 1$ (the extreme case $k=0$
is proved similarly). Then,  the matrix $q(\A(\G))$ has
maximum eigenvalue $q(\theta_0)>\eta_0$, with multiplicity  $1$.  Hence
the matrix
$\K:=\A(H)\otimes q(\A(\G))$ has eigenvalues
\begin{equation}
\label{ev-K}
\{q(\theta_0)\eta_e<q(\theta_0)\eta_{e-1}<\cdots<q(\theta_0)\eta_{e-i}
\le \eta_e\eta_0< \cdots < q(\theta_0)\eta_0\}
\end{equation}
for some $0\le i<e$; whereas  its weight-quotient
matrix $\B:=\A(H)\otimes \B_q$ ---with $\B_q$  being  the weight-quotient matrix of
$q(\A(\G))$ with respect to any partition $\{u\}\cup (V\setminus u)$---   has eigenvalues
\begin{equation}
\label{ev-Bq}
\textstyle
\{q(\theta_0)\eta_e< q(\theta_0)\eta_{e-1}<\cdots<q(\theta_0)\eta_{e-j}\le -\frac{\eta_0
w}{\|\vecnu\|^2-w}q(\theta_0)< \cdots< q(\theta_0)\eta_0\}
\end{equation}
for some $0\le j<e$. Furthermore, note that the multiplicity of $q(\theta_0)\eta_h$, $0\le
h\le e$, in (\ref{ev-K}) and (\ref{ev-Bq}) is the same ---and coincides with the multiplicity of
$\eta_h$ in $\A(H)$. Then, assuming that
$j\ge i$,   the eigenvalue inequality (coming again from Lemma \ref{interlac-coro})
$$
\eta_e\eta_0 \le - \frac{\eta_0 w}{\|\vecnu\|^2-w}q(\theta_0)
$$
gives the result. Now, it only remains to show that the other case is impossible.
Indeed,  if $j<i$, the same lemma would give
$$
q(\theta_0)\eta_{e-(j+1)} \le -\frac{\eta_0 w}{\|\vecnu\|^2-w}q(\theta_0)
$$
contradicting (\ref{ev-Bq}).
\endproof

From this result we can derive a number of consequences. For instance, notice that if we take
$H=K_r$, with $r=\alpha_k^w>1$, then
$1-\frac{\eta_0}{\eta_e}=r$ and (\ref{ineq-subgraphs}) becomes the
bound (\ref{(w,k)-indepen}) for the $(w,k)$-independence number.
As  other examples, we next give
some results, for regular graphs and distance-regular graphs, considering  the extreme cases
$k=0,1$ and $k=d-1$, respectively. First, if $\G$ is regular and $k=0$, we can take $H=\G$ and
Theorem  \ref{theo-ineq-subgraphs}
%(\ref{ineq-subgraphs})
---with $w=1$ and  $P_0(\theta_0)=1$--- gives again $1-\frac{\theta_0}{\theta_d} \le n$.
Still in the regular case, but taking now $k=1$, we also have the following corollary.

\begin{coro}
Let $\G$ be a regular  graph with $n$ vertices and
$\ev\G=\{\theta_0>\theta_1>\cdots>\theta_d\}$, $d\ge 2$, such that both $\G$ and its
complement $\overline{\G}$ are connected. Then,
$$
\frac{\theta_0-\theta_d}{\theta_1-\theta_d}\le
\frac{n}{n-(\theta_0-\theta_1)}\min\{-\theta_d,\theta_1+1\}.
$$
\end{coro}

\proof
Let $H=\G_{>1}\equiv \overline{\G}$. Then $\ev
H=\{\eta_0>\eta_1>\cdots >\eta_d\}$ with $\eta_0=n-\theta_0-1$ and
$\eta_i=-\theta_{d-i+1}-1$, $1\le i\le d$. Thus, using  (\ref{ineq-subgraphs}) with $w=1$ and
the value of $P_1(\theta_0)$ given in (\ref{P(theta)}), we get
$1+\frac{\theta_0-\theta_1}{\theta_1-\theta_d}
=\frac{\theta_0-\theta_d}{\theta_1-\theta_d}\le
\frac{n(\theta_1+1)}{n-(\theta_0-\theta_1)}$. Similarly, interchanging the roles of $\G$ and
$\overline{\G}$, we obtain $\frac{\theta_0-\theta_d}{\theta_1-\theta_d}\le
\frac{n(-\theta_d)}{n-(\theta_0-\theta_1)}$,  and the result follows.
%%(Of course, $\overline{\G}$ need not be connected, but it is certainly regular with the positive
%%eigenvector $\j$, so that Theorem \ref{theo-ineq-subgraphs} still holds.)
\endproof
%\newpage

Let us assume now that $\G$ is a distance-regular graph (see Biggs \cite{b93} or Brouwer et.
al. \cite{bcn89}). In this case, Theorem  \ref{theo-ineq-subgraphs} can be used to derive  the
following  upper bound for  the multiplicity of some of its eigenvalues.

\begin{coro}
\label{coro-mult-d-r}
Let $\G$ be a
%non-complete
distance-regular graph on $n$ vertices and with eigenvalues
$\ev \G=\{\theta_0>\theta_1>\cdots >\theta_d\}$. Then the multiplicity of some
eigenvalue $\theta_i$ with odd index $i$  satisfies the bound
\begin{equation}
\label{mult-d-r}
m(\theta_i)\le \frac{\pi_0}{\pi_i}\left( \frac{2n}{\sum_{j=0}^d\frac{\pi_0}{\pi_j}}-1\right)
\end{equation}
where $\pi_j:=\prod_{k=0,k\neq j}|\theta_j-\theta_k|$.
\end{coro}

\proof
Apply Theorem \ref{theo-ineq-subgraphs} with $H=\G_d$ ($k=d-1$) and $w=1$. Then, if $p_d$
denotes the distance-$d$ polynomial, satisfying $p_d(\A)=\A(\G_d)$, we have
$\eta_0=p_d(\theta_0)$ and $\eta_e={p_d}_{\min}:=\min_{1\le i\le d}p_d(\theta_i)$. Moreover,
by  (\ref{P(theta)}),  $P_{d-1}(\theta_0)+1=\sum_{j=0}^d\frac{\pi_0}{\pi_j}$. Consequently,
(\ref{ineq-subgraphs}) becomes
\begin{equation}
\label{ineq-d-r-graph}
\frac{p_d(\theta_0)}{-{p_d}_{\min}}\le \frac{2n}{\sum_{j=0}^d\frac{\pi_0}{\pi_j}}-1.
\end{equation}
Then, the result follows from  the known formula for the
multiplicities of a distance-regular graph in terms of $p_d$ (see, for instance, Bannai and Ito
\cite{bi84}),  namely
\begin{equation}
\label{mult-d-r-graph}
m(\theta_i)= (-1)^i\frac{\pi_0 p_d(\theta_0)}{\pi_i p_d(\theta_i)} \ \ \ (1\le i\le d),
\end{equation}
and the fact that, for some odd $i$, $1\le i\le d$,  we must have ${p_d}_{\min}=p_d(\theta_i)$ .
\endproof

This result is best possible in the sense that, for some distance-regular graphs, some of
the ``odd multiplicities" equal the upper bound in (\ref{mult-d-r}). In fact, we have examples
 where  all of such multiplicities equal the bound.
Indeed, if $\G$ is an $r$-antipodal distance-regular graph (see Biggs \cite{b93}), it was shown in
\cite{fgy4,g97} that
$P(\theta_0)+1=\sum_{j=0}^d\frac{\pi_0}{\pi_j}=2n/r$. Then, Corollary \ref{coro-mult-d-r}
assures that, for some odd index $i$, $1\le i\le d$, we have
$m(\theta_i)\le (r-1)\frac{\pi_0}{\pi_i}$ but, using that $\ev \G_d=\{r-1>-1\}$, it is easy to
prove that, in fact, $m(\theta_i)= (r-1)\frac{\pi_0}{\pi_i}$ for every $i=1,3,\ldots$ (see
\cite{f3,g97} for more details). In fact, in these references it was proved that a
distance-regular graph  is $r$-antipodal, for some $r\ge 2$, if and only if the multiplicities of
its eigenvalues are:
\begin{equation}
\label{mult-r-antipod-graph}
m (\theta_i)=\frac{\pi_0}{\pi_i} \quad \mbox{($i$ even)}, \ \ \ \ \ \
m (\theta_i)=(r-1)\frac{\pi_0}{\pi_i} \quad \mbox{($i$ odd)}.
\end{equation}

Notice that, since the distance-$d$ matrix $\A(\G_d)$ has positive eigenvector $\j$, with
eigenvalue $p_d(\theta_0)$, then $p_d(\theta_0)\ge |p_d(\theta_i)|$ for any $1\le i\le d$.
Consequently, (\ref{mult-d-r-graph}) gives the following general lower bound for the  ``even
multiplicities"
\begin{equation}
\label{mult-r-antipod-graph}
m (\theta_i)\ge \frac{\pi_0}{\pi_i}\quad \mbox{($i$ even)}
\end{equation}
 and the above example shows that this is also best possible. Also,
a direct proof of Corollary \ref{coro-mult-d-r} can be
obtained from (\ref{mult-r-antipod-graph}).

%%%%%%%%%%%%%%%%%%Sec. 7
\section{The Laplacian matrix}

When we deal with a non-regular graph $\G$, but still want to consider the cardinalities of
the vertex subsets, rather than their weights, we can use the Laplacian matrix $\L$ of
$\G$. As commented in Section 2, this is because it always has the eigenvalue $0$ with
eigenvector $\j$ (the
multiplicity of $0$ being the number of connected components of $\G$). Notice that
$\L$ can be seen as the adjacency matrix of a weighted pseudograph, obtained from $\G$ by
giving weight
$-1$ to its edges and adding a loop with weight $\delta_i$ on each vertex $v_i$. Therefore, as
when using the adjacency matrix $\A$, if the (distinct) vertices $u,v$ are $k$-independent, then
$(p(\L))_{uv}= 0$ for any polynomial $p$ of degree $k$. This allows us to derive some
results which are similar to those in the previous sections. For instance we next consider the
analogues of the bounds given in Sections 3 and 5.

\subsection{The independence number}

Let $\G$ be a $\delta$-regular graph on $n$ vertices, with adjacency matrix eigenvalues
$\lambda'_1=\delta\ge\lambda'_2\ge\cdots\ge\lambda'_n$. Since the Laplacian matrix of
$\G$ is $\L(\G)=\delta\I-\A(\G)$, its Laplacian eigenvalues are $\lambda_i=\delta-\lambda'_i$,
$1\le i\le n$, and hence the Hoffman-Lov\'asz' bound (\ref{alpha}) becomes
\begin{equation}
\alpha \le n\left(1-\frac{\delta}{\lambda_n}\right).
\end{equation}
In \cite{m88}  Mohar extended this bound to the case of non-regular graphs by considering
the degree average introduced below. In fact, as it is shown in the first part of the next
theorem,  Mohar's result can also be obtained reasoning as in the proof of Theorem
\ref{theo-weight-indepen}.
 Let $\delta_1\le \delta_2\le\cdots\le\delta_n$ represent the degree
sequence of $\G$ and  set $\db_r:=\frac{1}{r}\sum_{i=1}^r\delta_i$, $1\le r\le n$.

\begin{theo}
\label{theo-indepen-L}
Let $\G$ be a graph on $n$ vertices, with Laplacian eigenvalues
$\lambda_1=0\le\lambda_2\le\cdots
\le\lambda_n$. Then,
\begin{equation}
\label{indepen-L}
\alpha \le n\left(1-\frac{\db_{ \alpha}}{\lambda_n}\right).
\end{equation}
If the bound is attained for some independent set $C$, then $\G$ is $\delta$-regular, with
$\delta=\lambda_n(n-\alpha)/n$,  and
$C$ is a completely regular code, with every vertex in $\Cb$ being adjacent to
$\delta\alpha/(n-\alpha)=\lambda_n\alpha/n$ vertices of $C$.
\end{theo}

\proof
Let $C\subset V$ be a maximum independent set,  $\alpha=| C|$, with average degree
$\db_C:=\frac{1}{\alpha}\sigma_C=\frac{1}{\alpha}\sum_{u\in C}\delta_u$. Then,  the
weight-quotient matrix of $\L$ with respect to the partition $\Part$: $V_1\cup V_2=C\cup \Cb$
is now
\begin{equation}
\label{B-L}
\B=\left(\begin{array}{cc}
\frac{1}{\alpha}\sigma_C & \frac{-\sigma_C} {\sqrt{\alpha(n-\alpha)}} \\
\frac{-\sigma_C} {\sqrt{\alpha(n-\alpha)}} &
\frac{\sigma_C}{n-\alpha}
\end{array}\right)
=\db_C\left(\begin{array}{cc}
1 & \frac{-\sqrt{\alpha}} {\sqrt{n-\alpha}} \\
\frac{-\sqrt{\alpha}} {\sqrt{n-\alpha}} & \frac{\alpha}{n-\alpha}
\end{array}\right)
\end{equation}
with eigenvalues $\mu_1=0$  and  $\mu_2
%%= \tr \B =\frac{n}{\alpha(n-\alpha)}\sigma_C
=\frac{\db_C n}{n-\alpha}\le \lambda_n$,
by Lemma \ref{interlac-coro}, whence  (\ref{indepen-L}) follows since $\db_{\alpha}\le
\db_C$.

When equality holds $\db_\alpha=\db_C=\lambda_n(n-\alpha)/n$, the interlacing is
tight and, by Lemma \ref{interlac-coro}, the partition is pseudo-regular
with pseudo-quotient matrix (of $\L$ with respect to $\Part$)
$$
\B^*=\D\B\D^{-1}
=\delta\left(\begin{array}{cc}
1 & -1 \\
-\frac{\alpha} {n-\alpha} & \frac{\alpha}{n-\alpha}
\end{array}\right)
$$
where $\delta:=\db_\alpha$ and $\D=\diag(1/\sqrt{\alpha}, 1/\sqrt{n-\alpha})$.
But now the pseudo-intersection numbers of (\ref{intersec-num}) are simply
\begin{equation}
\label{intersec-num-L}
b_{ij}^*(u)=\left\{
\begin{array}{ll}\delta_u-\beta_{ij}(u), \mbox{ if } i= j \\
\ -\beta_{ij}(u), \mbox{ otherwise }
\end{array}
\right.
\ \ \ \ (u\in V_i)
\end{equation}
%(recall that each vertex $u$ is adjacent to itself with weight $\delta_u$)
where $\beta_{ij}(u):=|\G(u)\cap V_j|$
represent the (standard) intersections numbers of $\G$ with respect to $\Part$. This gives:
$\beta_{11}(u)=0=\delta_u-\delta$,   whence
$\delta_u=\delta$,
$\beta_{12}(u)=\delta$ for any $u\in C$; and $\beta_{21}(v)=\frac{\delta\alpha}{n-\alpha}$,
$\beta_{22}(v)=\delta-\frac{\delta\alpha}{n-\alpha}$ for any
$v\in \Cb$; whence the result follows.
\endproof

In order to make the inequality (\ref{indepen-L}) more explicit, Mohar \cite{m88} presents his
result by stating that, if $r$ is the smallest positive integer for which $r>
n(\lambda_n-\db_r)/\lambda_n$, then $\alpha \le r-1$.

The second part of the theorem, characterizing the case of equality, extends a result of Haemers
\cite{ha95} for regular graphs. (He uses the adjacency matrix and hence considers the
---equivalent--- case of equality in the Hoffman-Lov\'asz' bound (\ref{alpha}) ---see the last
comment in the proof of Theorem \ref{theo-weight-indepen}.)

\subsection{The set independence number}
Let us now see how the results of Section 5 look when the Laplacian spectrum is involved.
We shall here omit the proofs, since they are very similar, and the use of the Laplacian matrix
has already been illustrated above. Of course, since
$\ev\L=\{\theta_0=0<\theta_1< \cdots<\theta_d\}$, the
$k$-alternating  polynomial $P_k$, defined as in (\ref{def-P_k}), must now attain maximum
value at $0$, that is on the left of the other eigenvalues.

\begin{theo}
\label{theo-(w,k)-indepen-L}
Let $\G$ be a connected graph on $n$ vertices, with Laplacian eigenvalues  $\ev \L=\{0<
\theta_1<\cdots < \theta_d\}$,  and corresponding $k$-alternating polynomials $P_k$, $0\le k\le
d-1$.  Assume that, for some  integer $w$, $1\le w\le n/2$, the
$(w,k)$-independence number satisfies $\alpha_k^w\ge 2$. Then,
\begin{equation}
\label{(w,k)-indepen-L}
\alpha_k^w \le  \frac{2n}{w(P_k(0)+1)}.
\end{equation}
\end{theo}

In particular, using  that
$P_1(0)=2\frac{-\theta_1}{\theta_1-\theta_d}+1
=\frac{\theta_d+\theta_1}{\theta_d-\theta_1}$,
we have that the independence and chromatic numbers
%$\alpha\equiv \alpha_1^1$
of a connected graph $\G$, in terms of its Laplacian eigenvalues, satisfy respectively
$$
\alpha \le n\left(1-\frac{\theta_1}{\theta_d}\right), \ \ \ \ \ \
\chi\ge \frac{\theta_d}{\theta_d-\theta_1}.
$$
Of course, the advantage of the first bound above, in comparison with (\ref{indepen-L}), is its
%exclusive  dependence on the Laplacian spectrum.
explicit form in terms of the Laplacian spectrum.
\begin{coro}
Let $\G$ be a connected graph on $n$ vertices, $(^n_w)\ge t$, with Laplacian eigenvalues
$\ev \L=\{0<\theta_1<\cdots<\theta_d\}$. Then,
\begin{equation}
\label{coro-wt-diam-L}
P_k (0) > \frac{2n}{wt}-1 \ \  \Rightarrow\ \  \ D_t^w\le k .
\end{equation}
\end{coro}
As in Section 5,  in the above results we can replace $P_k$ by the Chebychev polynomial
$T_k(-x)$  ``shifted" from $[-1,1]$ to $[\theta_1,\theta_d]$, that is
$T_k\left( \frac{\theta_1+\theta_d-2x}{\theta_d-\theta_1} \right)$, now giving:
\begin{equation}
\label{Dtw-L}
D_t^w \leq \left\lfloor\frac{{\rm cosh}^{-1} \left(\frac{2n}{wt}-1\right)}
 {{\rm cosh}^{-1}\left(\frac{\theta_d
+ \theta_1}{\theta_d-\theta_1} \right)} \right\rfloor +1,
\end{equation}
whereas the results of Chung, Delorme, and Sol\'e \cite{cds97},  proved by using the normalized
Laplacian  matrix (the so-called {\it Laplace operator}), correspond to
\begin{equation}
\label{Dtw-L-chung}
D_t^w \leq \left\lfloor\frac{{\rm cosh}^{-1}
\left(\frac{n}{w}-1\right)} {{\rm cosh}^{-1}\left(\frac{\theta_d
+\theta_s}{\theta_d-\theta_s} \right)} \right\rfloor +1,
\end{equation}
where $\theta_s$ is the $t$-th smallest Laplacian
eigenvalue $\lambda_t$, that is, $s$ is the smallest integer satisfying
$1+m(\theta_1)+\cdots+m(\theta_s)\ge t$.
%in $\lambda_1=0<\lambda_2\le\cdots\le \lambda _n$
%are the $n$ eigenvalues of the Laplacian matrix $\L$
See also Chung, Grigor'yan, and Yau
\cite{cgy96}. In the way of comparing the above bounds, note that for $t=2$ both results
coincide. Otherwise, there is a general case in which  (\ref{Dtw-L}) clearly supersedes
(\ref{Dtw-L-chung}), namely whenever the multiplicity of the second smallest Laplacian
eigenvalue $\theta_1$ ---the so-called ``algebraic connectivity" of $\G$--- satisfies
$m(\theta_1)\ge t-1$ (since then $s=1$).
%$\lambda_t=\theta_1$).
%%%%%%%%%%%%%%%%%%%%%%%
\vskip .4cm

\noindent{\bf Acknowledgment.}
%Work supported in part by the Spanish Research Council
%(Comi\-si\'on Interministerial de Ciencia y Tecnolog\'\i a, CICYT)
%under projects TIC 92-1228-E and TIC 94-0592.
I am indebted to the referee for helpful comments and suggestions which lead to numerous
improvements of the manuscript.

%\newpage
%%%%%%%%%%%%%%%%%Ref.


\begin{thebibliography}{99}

\bibitem{bcff95}
C. Balbuena, A. Carmona, J. F\`abrega, and M.A. Fiol, On the connectivity and the
conditional  diameter of graphs and digraphs, {\it  Networks} {\bf 28} (1996), 97--105.

\bibitem{bi84}
E. Bannai and T. Ito, ``Algebraic Combinatorics I: Association Schemes," Benjamin-Cummings
Lecture Note Ser. {\bf 58}, Benjamin/Cummings, London, 1984.

\bibitem{b79}
N  Biggs, Some odd graph theory,  {\it Annals New York Acad. Sci.} {\bf 319}, New York,
1979,  71--81.

\bibitem{b93}
N. Biggs, ``Algebraic Graph Theory," Cambridge Univ. Press, Cambridge, 1974; 2nd ed., 1993.

\bibitem{bd88}
I. Bond and C. Delorme,  New large bipartite graphs with given degree and diameter,
{\it Ars Combin.} {\bf 25C} (1988), 123--132.

\bibitem{bcn89}
 A.E. Brouwer, A.M. Cohen, and A. Neumaier, ``Distance-Regular Graphs,"
Springer-Verlag, Berlin, 1989.

\bibitem{cds97}
F.R.K. Chung, C. Delorme, and P. Sol\'e, $k$-Diameter and spectral multiplicity,
submitted.

\bibitem{cgy96}
F.R.K. Chung, A. Grigor'yan, and S.-T. Yau, Upper bounds for eigenvalues for the discrete and
continuous Laplace operators, {\it Adv. Math.} {\bf 117} (1996), 165--178.

\bibitem{CH}
R. Courant and D. Hilbert, ``Methoden der Mathematischen Physik,"  Vol. 1, Springer-Verlag,
Berlin, 1924.

\bibitem{vd96}
E.R. van Dam, ``Graphs with Few Eigenvalues," Ph.D. Thesis, Tilburg University, Tilburg, 1996.

\bibitem{vd98}
E.R. van Dam, Bounds on special subsets in graphs, eigenvalues and association schemes, {\it J.
Algebraic Combin.} {\bf 7} (1998),  no. 3,  321--332.

\bibitem{vdh95}
E.R. van Dam and W.H. Haemers, Eigenvalues and the diameter of graphs,
{\it Linear and Multilinear Algebra} {\bf 39} (1995), 33--44.

\bibitem{dt97}
C. Delorme and J.P. Tillich, Eigenvalues, eigenspaces and distances to subsets, {\it Discrete
Math.} {\bf 165-166} (1997), no. 1-3, 161--184.

\bibitem{f1}
M.A. Fiol, Weight odd parameters and spectra of graphs, {\it in} ``Eighth
Quadrennial International Conference on Graph Theory, Combinatorics, Algorithms,
and Applications," Kalamazoo, MI, USA, June 3-7, 1996.

\bibitem{f3}
M.A. Fiol, An eigenvalue characterization of antipodal distance-regular graphs, {\it Electron. J.
Combin.} {\bf 4} (1997),  no. 1, \#R30.

\bibitem{fg1}
M.A. Fiol and E.  Garriga, The alternating and adjacency polynomials,  and their  relation with the
spectra and diameters of graphs, {\it Discrete Appl. Math.} {\bf 87}
(1998), no. 1-3, 77--97.

\bibitem{fg3}
M.A. Fiol and E. Garriga,
On the algebraic theory of pseudo-distance-regularity around a set, {\em Linear Algebra Appl.} {\bf 298} (1999), no. 1-3, 115--141.

\bibitem{fgy1}
M.A. Fiol,  E. Garriga, and J.L.A. Yebra,  On a class of polynomials and its relation with the spectra and diameters of graphs,  {\it J.  Combin. Theory Ser. B} {\bf 67} (1996), 48--61.

\bibitem{fgy4}
M.A. Fiol, E. Garriga, and J.L.A. Yebra, Boundary graphs: The limit case of a spectral property, {\it Discrete Math.} {\bf 226} (2001), no.1-3, 155--173.

\bibitem{fgy5}
M.A. Fiol,  E. Garriga, and J.L.A. Yebra, The alternating polynomials and their relation with the spectra and conditional diameters of graphs, {\it Discrete Math.} {\bf 167--168} (1997), no. 1-3, 297--307.

\bibitem{g97}
E. Garriga, ``Contribuci\'o a la Teoria Espectral de Grafs. Problemes M\`etrics i Dist\`ancia-Regularitat,"  Ph.D. Thesis, Universitat Polit\`ecnica de Catalunya, Barcelona, 1997.

\bibitem{g93}
C.D. Godsil, ``Algebraic Combinatorics,"  Chapman and Hall, London/New York, 1993.

\bibitem{gm80}
C.D. Godsil and B.D. McKay, Feasibility conditions for the existence of walk-regular graphs,
{\it Linear Algebra Appl.} {\bf 30} (1980), 51--61.

\bibitem{ha79}
W.H. Haemers,  Eigenvalue methods, {\it in} ``Packing and Covering in Combinatorics" (A.
Schrijver, Ed.), Math. Centre Tract {\bf 106}, Mathematical Centre, Amsterdam, 1979, 15--38.

\bibitem{ha80}
W.H. Haemers, ``Eigenvalue Techniques in Design and Graph Theory,"  Math. Centre Tract
{\bf 121}, Mathematical Centre, Amsterdam, 1980.

\bibitem{ha95}
W.H. Haemers, Interlacing eigenvalues and graphs, {\it Linear Algebra Appl.}
{\bf 226-228} (1995), 593--616.

\bibitem{hof70}
A.J.  Hoffman,  On eigenvalues and coloring of graphs, {\it in}  ``Graph Theory and its
Applications," (B. Harris, Ed.), Academic Press, New York, 1970, 79--91.

\bibitem{k97}
N. Kahale, Isoperimetric inequalities and eigenvalues, {\it Siam J. Discrete Math.} {\bf 10}
(1997), no. 1, 30--40.

\bibitem{knu94}
D.K. Knuth, The sandwich theorem,  {\it Electron. J. Combin.} {\bf 1} (1994), \#A1.

\bibitem{l79}
L. Lov\'asz,  On the Shannon capacity of a graph, {\it IEEE Trans. Inform. Theory}
{\bf IT-25} (1979), no. 1, 1--7.

\bibitem{l79-2}
L. Lov\'asz,  ``Combinatorial Problems and Exercises," North-Holland, Amsterdam, 1979.

\bibitem{m76}
B.D. McKay, ``Backtrack Programming and the Graph Isomorphism Problem,"  M.Sc. Thesis,
University of Melbourne, 1976.

\bibitem{m88}
 B. Mohar, The Laplacian spectrum of graphs, {\it in} ``Graph Theory,
Combinatorics and Applications" (Y. Alavi et.al., Ed.) Vol. 2, John Wiley \& Sons, New York, 1991,
871--898.

\bibitem{r97}
J.A. Rodr\'\i guez, ``Cotas de Diversos Par\'ametros de un Grafo a Partir de los Autovalores de su
Matriz Laplaciana"  Ph.D. Thesis, Universitat Polit\`ecnica de Catalunya, Barcelona, 1997.

\bibitem{sha56}
C.E. Shannon, The zero-error capacity of a noisy channel, {\it IRE Trans. Inform. Theory} {\bf
IT-2} (1956), no. 3, 8--19.

\end{thebibliography}
\end{document}